# Optimal queueing-based rebalancing for one-way electric carsharing systems with stochastic demand

Tai-Yu Ma[1], Theodoros Pantelidis[2], Joseph Y. J. Chow[2]

[1]Luxembourg Institute of Socio-Economic Research
11 Porte des Sciences, L-4366 Esch-sur-Alzette, Luxembourg

[2]Department of Civil & Urban Engineering, Tandon School of Engineering, New York University
6 Metro Tech Center, Brooklyn, NY, 11201, USA

## ABSTRACT

Viability of electric vehicle car sharing operations depends on rebalancing algorithms. Earlier methods in the literature suggest a trend toward Markovian stochastic demand with server relocation with queueing constraints. We propose a new model formulation based on a node-charge graph structure that extends the relocation model to include transshipment relocation flows. Computational tests with up to 1000 node (and 4000 node-charges) suggest promising avenues for further study.

# 1. INTRODUCTION

Car sharing operations are an essential part of "smart mobility" solutions in congested large megacities. According to Martin and Shaheen (*1*) a single carshare vehicle can replace 7 to 11 personal vehicles on the road. The common practice in such services is to book specific time slots and reserve a vehicle from a specific location. The return location is required to be the same for "two-way" systems but is relaxed for "one-way" systems. Within one-way systems, "station-based" systems restrict vehicles to specific parking locations while "free-floating" systems allow returns anywhere within a covered area. Examples of free floating systems are the BMW ReachNow car sharing system in Brooklyn and Car2Go in New York City with service areas in 2017 shown in Figure 1 (the two companies merged in 2018 (*2*)).

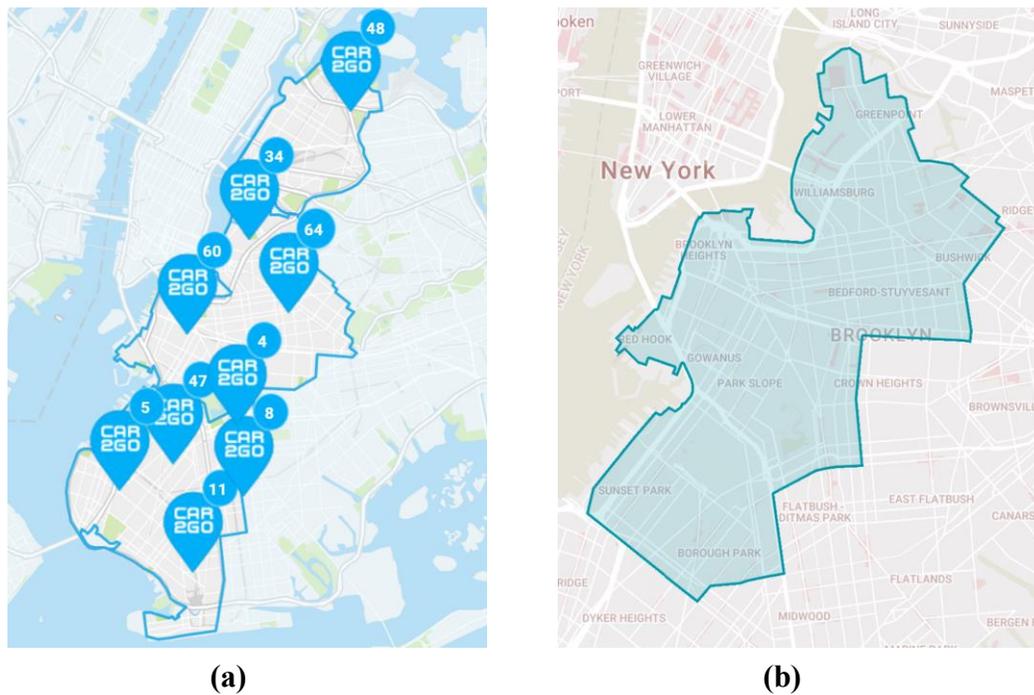

(a) (b)
**FIGURE 1.** Examples of free-floating carshare systems: (a) Car2Go (source: car2go.com) and (b) BMW ReachNow (source: reachnow.com).

In large car sharing systems, vehicle rebalancing is one of the primary challenges to ensure efficiency and an adequate level of service. Potential customers may end up waiting or accessing a farther location, or even balk from using the service, if there is no available vehicle within proximity, or no parking or return location available near the destination. Rebalancing involves having either the system staff or users (through incentives) drop off vehicles at locations that would better match supply to demand.

Many studies have been devoted to the rebalancing problem (*3 – 7*) and related regulating challenges: station location with rebalancing (*8*); pricing incentives (*9 – 11*); parking reservations (*12*); routing personnel (*6, 13*); parking capacity, fleet sizing (*14*), etc.

Carsharing companies have further considered electric vehicle (EV) fleets to be more sustainable and to reduce gasoline consumption costs. Not only do Zipcar, Car2Go, and ReachNow all operate some EVs in their worldwide fleets, some startup carsharing businesses rely exclusively on EVs: e.g. Autolib' and Cité Lib in France, BlueIndy, DriveNow in Copenhagen, Carma based in San Francisco, and Los Angeles' Low-Income Plan for EV car sharing (*15*). However, rebalancing is further hampered in an EV environment. Under such a setting, vehicles not only need to be rebalanced in space to meet stochastic customer demand, but they also need to be rebalanced over time in terms of charging to meet future customer demand. A smaller subset of studies emerged in recent years to tackle this heightened challenge. Two general methods have been adopted. The first assumes demand is sufficiently deterministic in a multiperiod setting (e.g. (*16*)). This can be problematic for systems in which most of the demand is not made for repeated commute trips and/or the fleet is spatially distributed in such a manner that the vehicle density is fairly low, which is most systems. The second group of methods assumes stochastic demand, either through stochastic programming (*17*), simulation (*18*), or with Markovian demand (*19*). The latter Markovian queueing models appear promising, but earlier EV studies either assume a simplistic relocation policy (*19*) or, in the case of queueing networks (*10*), do not allow customers to pick up vehicles at nearby locations. Discrete network approaches have not been used with queueing models for EV charging setting either.

It becomes quite clear that EV car sharing systems present a more complex environment in terms of mathematical modeling and decision analysis. The main goal of our study is to propose a new mathematical model for rebalancing carsharing vehicles over time that: 1) captures Markovian queueing in customer demand and 2) allows customers to pick up vehicles at different charge levels and at different locations further than their own "zone". The model is structured as a p-median facility relocation problem with queueing constraints similar to Marianov and Serra (*20*) and Sayarshad and Chow (*7*), but with the extension of the relocation assignment component from a bipartite transportation problem into a single-commodity minimum cost flow (transshipment) assignment problem.

The rest of the paper is organized as follows: Section 2 presents the literature review of previous studies in the field; Section 3 introduces the modeling framework and a small demonstration of the MIP formulation that is proposed for large networks as shown in Section 4. Section 5 provides an evaluation of these results and main research findings; Section 6 includes concluding remarks and suggestions for future research.

## 2. LITERATURE REVIEW

Early studies for designing carsharing systems relied on simulation (*21*) for evaluation. More systematic mathematical models to optimizing carsharing fleets have since been proposed with rebalancing in mind. The purpose of optimal rebalancing is to make the decision considering trade-offs among a set of parameters that differ from one instance to another:

- Customer demand for vehicles over time, with distribution for pickup and drop-off locations
- Vehicle access distance for customers, which should not be restricted to a specific "zone"
- Penalty for having insufficient vehicle stock to meet demand
- Cost of relocation
- Distribution of booking duration, which may exceed one direct trip from pickup to drop-off – for example, a customer may take a vehicle out of town to run errands before returning the vehicle at another location in the coverage area
- Rental cost over time
- Distribution of remaining fuel upon return of vehicle

- For station-based one-way systems, capacity of spaces in each zone is a factor
- Reservation time requirement and availability of information to the system

In addition, for EV systems there are additional trade-offs:
- Cost of charging EVs
- Duration of EV charging
- Demand for EVs at a minimum charge level
- Distribution of charge consumed upon return of vehicle
- Capacity of charging stations, i.e. number of chargers

The evolution of carsharing rebalancing models can be broken down into four general stages with increasing ability to address key trade-offs: deterministic demand assignment, stochastic demand, queueing network, and queueing-based facility location. Simulation methods are left out as we are interested in analytical solutions. We summarize some example representative studies (non-EV and EV) as well as the trade-offs that they handle well or poorly. Note that they do not necessarily follow a chronological order as different methods from different stages have been used for different applications.

**Deterministic demand assignment**: these are characterized by network models, sometimes with multiple time periods to capture capacity costs in terms of delay costs for customers.
- Non-EV studies: Febbraro et al. (*22*), Chow and Sayarshad (*8*), Jorge et al. (*9*), Kaspi et al. (*12*)
- EV studies: Boyaci et al. (*23*), Xu et al. (*16*)
- Weaknesses: does not handle random demand or its consequences like stockout/balking penalties, information availability to the system, and impacts of charging station capacities

**Stochastic demand**: these models consider stochastic representations of the demand, in some cases involving Markovian demand with one-step lookahead via two-stage stochastic programming or chance constraints
- Non-EV studies: Nair and Miller-Hooks (*5*)
- EV studies: Brandstätter et al. (*17*)
- Weaknesses: while demand is random, the decisions are generally myopic (ignores dependency of future decisions on current decisions) and reliant on simulation-based optimization for capturing all the trade-offs

**Queueing network**: these models consider the steady state impacts of rebalancing decisions under a stochastic environment, which are more dynamically stable
- Non-EV studies: Waserhole and Jost (*10*), Zhang and Pavone (*24*)
- EV studies: none published yet
- Weaknesses: the restriction to have demand be served only at the specified origin/destination zones assumes customers are not able to switch to a different location nearby, which is not realistic for one-way systems, particularly free-floating systems

**Queueing-based facility location**: these are characterized by facility location models that incorporate queueing criteria in the objectives or constraints, and allow for interzonal access costs

Ma, Pantelidis, Chow 5- Non-EV studies: Sayarshad and Chow (*7*), Ma et al. (*25*)
- EV studies: none published yet
- Weaknesses: see below

Facility location can consider queue delay by treaching each service node as a queue with the the number of servers $s$. In this case, however, even a simple assumption of an M/M/s stochastic queue results in a nonlinear objective. Since nonlinear integer programming problems are undesirable, researchers have proposed alternative methods to handle the queueing. One such ways is the Q-MALP model from Marianov and ReVelle (*26*), who showed that the queue delay objective can instead be cast as a set of piecewise linear constraints for the intensity to be within a specified reliability level $\eta$. Because the intensity parameter can be preprocessed for different numbers of servers, it is possible to solve a facility location with desired queueing-based service reliability as a mixed integer linear programming problem. The model has since been modified to handle maximal coverage (*27*), server allocation (*20*), and p-median coverage with relocation costs (*7*).

Queueing-based facility location models handle everything that the "Queueing network" models can, and on top of that they allow for interzonal matching of vehicles to demand. However, the relocation component is based on a bipartite transportation problem of moving excess servers to locations in demand of servers. This is fine for a non-EV carsharing system, but for EV charging the mechanics are more complex (see (*28*)). The model does not distinguish demand for a minimum charge level that can be less than 100%. For example, a customer should be allowed to demand a vehicle with 60% charge, and this same customer should be allowed to pickup a vehicle with 80% charge.

We address this issue by proposing the first queueing-based facility relocation problem and formulate it as a mixed integer programming problem that can be solved sufficiently fast even for 100-zone examples using commercial integer programming solvers. The model hinges on defining a new graph structure.

## 3. PROPOSED MODEL

### 3.1 Node-charge graph structure

Before diving into the model formulation of the EV rebalancing problem, we need to modify the conventional graph structure such that the trade-offs can all be addressed. A network of zones is connected by links. A subset of these zones is designated as a set of charging stations with finite numbers of chargers. At the start of every time interval, there is a set of idle vehicles that have not been reserved by any customer. These vehicles may either be sitting somewhere in a zone unused, charging at a charging station, or being relocated to another zone or charging station. The locations and charge levels of the vehicles are known. Customers arrive randomly within that time interval according to a stationary Poisson process. When they book the vehicles for use, the vehicles are effectively "serving" the customers for a period of time that is assumed to follow an exponential distribution. Customer arrivals and vehicle return zone locations are assumed to follow discrete distributions.

We can graphically illustrate this in a one-dimensional network without loss of generality. Consider a 5-node network lined up in sequence as shown in Figure 2(a) where node 1 and node 3 are charging stations (denoted by gray nodes).

The graph is converted into a node-charge graph representation in Figure 2(b) where the y-axis is a discrete charge interval (let's say there are four intervals: 20%+, 40%+, 60%+, 80%+). Each layer represents the same zones at a certain charge level. Unidirectional links exist at the charging stations to represent charging with link costs representing charging cost and time. The charging links are also capacitated. A vehicle positioned at a node covers ALL nodes underneath it with lower charge as well; access costs for demand is based only on the spatial link costs and not the charging costs. For example, a vehicle at node 4 with charge 40+ can serve node 1 … 5 at 40%+ and also at 20%+, as illustrated by the blue arrows. The access cost of demand at (node 2, charge 20+) for the vehicle at (node 4, charge 40+) is just the cost from node 2 to node 4.

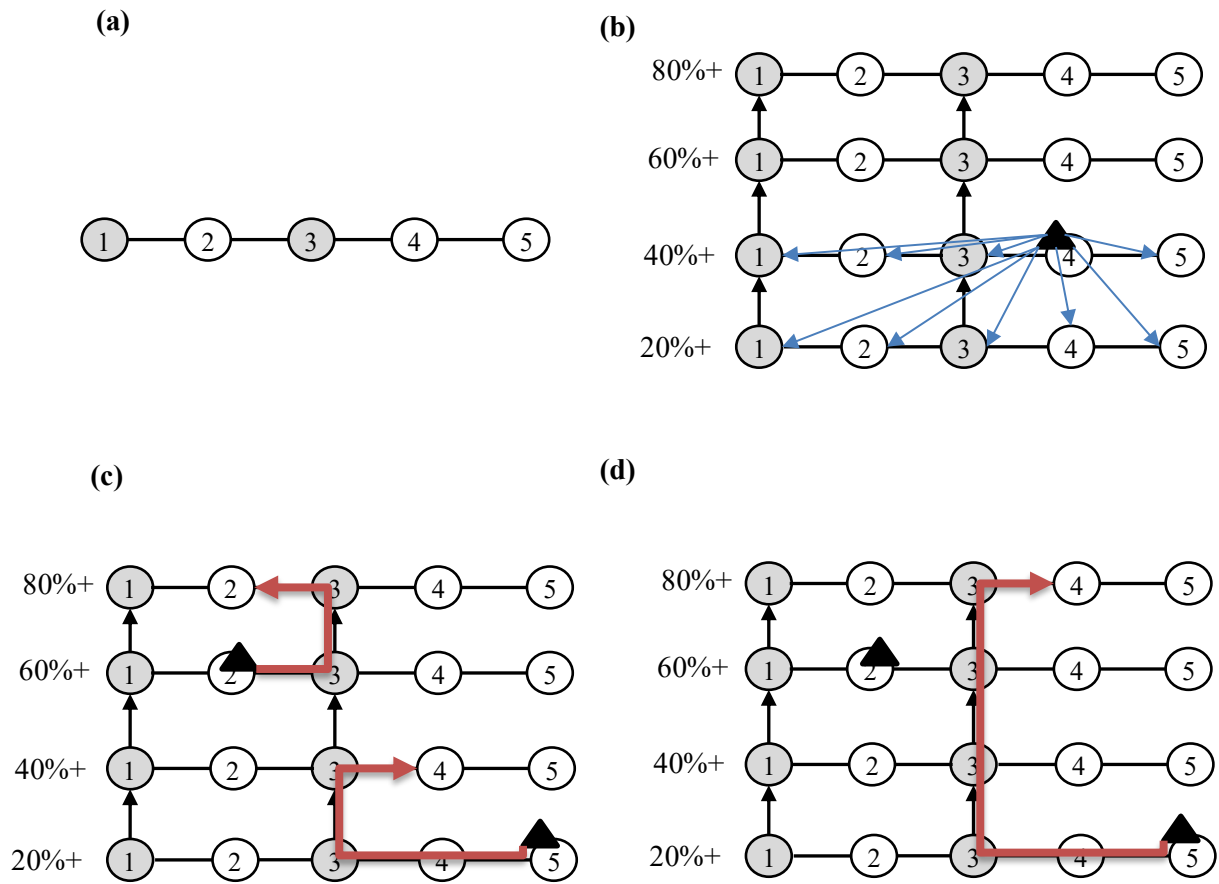

**FIGURE 2.** (a) initial graph; (b) expansion to a node-charge graph with coverage illustration in blue arrows; (c) one possible rebalancing solution in red arrow; and (d) a second rebalancing solution.

Suppose there are two idle vehicles, one with 20% charge at node 5 and one with 60% charge at node 2. Two possible optimal rebalancing solutions are shown in Figure 2(c) and 2(d). In Figure 2(c), one vehicle is directed toward station 3 to recharge 20% before being relocated back to node 2, while the second vehicle is sent to station 3 to recharge 20% before being relocated to node 4. In Figure 2(d), only one vehicle is relocated to station 3 to be recharged up to 80%+ and would then be relocated to node 4.

The optimality of these solutions depends on a mix of factors. For example, Figure 2(c) may be best if there is high demand near nodes 1, 2, and 3 up to 80%+ charge and there is enough capacity at station 3 to allow two vehicles to charge at the same time. The charging cost might be very high relative to the spatial relocation and/or access costs/penalties, leading to two short charges instead of one longer charge. Figure 2(d) may be best if the high charge demand is located closer to nodes 4 and 5, and/or perhaps there is only enough capacity for 1 vehicle charging at station 3 and the relocation cost to station 1 does not warrant the additional charging.

This graph structure is highly flexible. For example, it can be set to 100 discrete charge levels. Duplicates of locations can be used (e.g. a zone 3 with a zone "3a" for charging stations at zone 3) to model the decisions to remove a charged vehicle from a station without having to rebalance it to another zone. The bottom charge level can represent minimum charge needed to get from any zone to any other zone without depleting fully.

Under this graph structure, we model the relocation decision such that demand at one or more charge levels are assigned to each zone, where the service at the zone acts as a stochastic queue with one or more servers.

**3.2 Model formulation**

The following notation is used.

$N$: set of nodes
$H$: set of charging levels, $H = 1,2,...$
$A$: set of directed arcs in the node-charge graph, $A = \{((i,g),(j,h)) | \forall i,j \in N, g,h \in H\}$
$B$: total number of idle vehicles at the start of a rebalancing time interval
$y_{ig}$: number of idle vehicles at the node-charge $(i,g)$ at the start of a rebalancing time interval
$J \subset N$: subset of nodes that are charging stations
$\theta$: a conversion parameter
$u_j, j \in J$: capacity of charging facilities available
$c_{igjh}$: cost on an arc of node-charge graph from $(i,g)$ to $(j,h)$
$O$: set of origins of idle vehicles on the node-charge graph
$\lambda_{ih}$: arrival rate of customers at node $i$ with demand for SOC $h$ intervals or higher; assume that customers use up exactly $h$ intervals during their trip;
Parameters should be defined such that $\lambda_{idh} = 0$ if $c_{id} > h\varphi$, where $\varphi$ is the amount of charge in one interval
$C$: max number of servers at a node-charge
M: large positive penalty constant
$A_{ig}^+$: set of outgoing arcs originating at a node-charge $(i,g), i \in N, g \in H$
$A_{ig}^-$: set of incoming arcs destinated to a node-charge $(i,g), i \in N, g \in H$
$W_{jgjh}$: rebalancing EV flow on arc $((i,g),(j,h)) \in A$,
$\boldsymbol{Y_{jhm}}$: if the m-th idle vehicle are located at node $j$ with charge $h$
$\boldsymbol{X_{igjh}}$: vehicle at node $j$ with charge level $h$ serves customer demanding $g \in \{1,...,|H|\}$ charge or more at node $i$ if $X_{igjh} = 1$
$\underline{p_{jgjh}}$: Partial path flow from node-charge $(j,g)$ to $(j,h), \forall j \in J, h > g,\ g,h \in H$

The objective of the problem is to minimize total customer's cost and the operator's weighted rebalancing cost under customer stochastic demand and charging station capacity constraints. The problem is an extension of a multiple server relocation problem under stochastic demand (7) in the context of electric vehicle charging. We consider the problem on a digraph with multilevels (multilayers) that converts the EV-rebalancing problem into a facility location embedded with a single commodity minimum cost flow problem to rebalance idle vehicles to meet customer demand for different charging levels.

Let $G(N, A)$ be a directed graph with $N$ being a set of node-charge $(i, h), \forall i \in N, h \in H$, and $A$ a set of directed arcs, $A = \{((i, g), (j, h)) | \forall i, j \in N, \forall g, h \in H\}$. A node-charge $(i, h)$ is characterized by demand's spatial location $i \in N$ and request charging level (battery level) $g \in H$. We discretize customer's charging demand into a set of charging levels $H$. All demand between two levels $h-1$ and $h$ sum up to level $h$. Arcs only cross from one charging level up to another at charging station nodes $J$ and there is a capacity $u_j > 0, j \in J$ applied to all flows going through that station node regardless of charge level. Charging arcs are a subset of arcs defined as $A_{\hbar} = \{((j, g), (j, h)) | \forall j \in J, \forall h = g+1, g = 1, \ldots |H-1|\}, A_{\hbar} \subset A$. The assigned flow on arcs is an integer decision variable. The cost of assigning flow on arcs is the multiplication of arc flows by its unitary cost, measured as rebalancing time/cost or charging time/cost on arcs.

Rebalancing idle vehicles is making tradeoffs between serving customer demand at higher access cost, increasing relocation cost to reduce access cost, and increasing charging cost to cover demand for higher charged vehicles at the expense of the vehicles' availability. The problem is formulated as a p-median facility location problem embedded with a capacitated commodity flow problem. The idle vehicle relocation problem in an electric carsharing system is formulated as follows.

$$\min Z = \sum_{i \in N} \sum_{j \in N} \sum_{h \in H} \sum_{g \in H} \lambda_{ig} t_{ij} X_{igjh} + \theta \sum_{((i,g),(j,h)) \in A} c_{igjh} W_{igjh} \quad (1)$$

s.t.

$$\sum_{j \in N} \sum_{h \in H, h \geq g} X_{igjh} = 1, \quad \forall i \in N, g \in H \quad (2)$$

$$\sum_{j \in N} \sum_{h \in H, h < g} X_{igjh} = 0, \quad \forall i \in N, g \in H \quad (3)$$

$$Y_{jhm} \leq Y_{jh,m-1}, \quad \forall j \in N, h \in H, m = 2, 3, \ldots, C \quad (4)$$

$$\sum_{i \in N} \sum_{g \in H} \lambda_{ig} X_{igjh} \leq \mu_{jh} \left[ Y_{jh1} \rho_{\eta jh1} + \sum_{m=2}^{C} Y_{jhm} (\rho_{\eta jhm} - \rho_{\eta jh,m-1}) \right], \forall j \in N, h \in H \quad (5)$$

$$\sum_{j \in N} \sum_{h \in H} \sum_{m=1}^{C} Y_{jhm} = B \quad (6)$$

$$X_{igjh} \leq Y_{jh1}, \quad \forall i, j \in N, g, h \in H \quad (7)$$

$$\sum_{(j,h) \in A_{ig}^-} W_{igjh} - \sum_{(j,h) \in A_{ig}^+} W_{igjh} \leq M Y_{ig1}, \forall (i, g) \in N \times H \backslash O \quad (8)$$

$$-(\sum_{(j,h)\in A_{ig}^-} W_{igjh} - \sum_{(j,h)\in A_{ig}^+} W_{igjh}) \leq MY_{ig1}, \forall (i,g) \in N \times H\setminus O \tag{9}$$

$$\sum_{(i,g)\in A_{jh}^-} W_{igjh} - \sum_{(j,h)\in A_{ig}^+} W_{igjh} + y_{ig} = \sum_{m=1}^{C} Y_{jhm}, \quad \forall j \in N, h \in H \tag{10}$$

$$\sum_{g'\leq g}\sum_{h'\geq h} p_{jg'jh'} = W_{jgjh} \quad \forall j \in J, \forall ((j,g),(j,h)) \in A \tag{11}$$

$$\sum_{g=1}^{|H|-1}\sum_{h'=g+1}^{H} p_{jgjh'} \leq u_j, \forall j \in J, \forall ((j,g),(j,h)) \in A \tag{12}$$

$$X_{igjh} \in \{0,1\}, \forall i,j \in N, g,h \in H \tag{13}$$

$$Y_{jhm} \in \{0,1\}, \forall j \in N, h \in H, m = 1,2,3,\dots,C \tag{14}$$

$$W_{igjh} \in 0 \cup Z^+, \forall i,j \in N, g,h \in H \tag{15}$$

$$p_{jgjh} \geq 0, \forall ((j,g),(j,h)) \in P_j \tag{16}$$

The objective function minimizes the total access cost of customers to servers (idle vehicles) and total routing cost of idle vehicles (i.e. travel time/cost from the current locations of idle vehicles to charging stations, charging time/cost and travel time/cost to its respective destinations) on the node-charge graph (network). The rebalancing operations are run at each predefined time interval in order to serve customer's demand and minimize queueing delay and operating cost.

Constraints (2) and (3) require that rebalanced idle vehicles to serve randomly arriving customers must have sufficient charge to match the demanded amount. Constraint (4) is an order constraint stating a m-th server can be present only if there is already a (m-1)-th server at the same location.

Constraint (5) is the piecewise linear queueing constraint from Marianov and ReVelle (*26*) queueing constraint representing when a customer arrives at an idle vehicle, there will be no more than $b$ other customers waiting on a line with a probability more than service reliability $\eta$. The intensity is setup as a recursive cumulative value based on the number of servers assigned to the location.

Constraint (6) states the total number of servers is equal to the total number of available idle vehicles. Constraint (7) assures that only a location with servers can cover demand nodes. Constraints (8-10) are the flow conservation constraints of the minimum cost flow problem. For the charging station capacity, we need to assure that the assigned flow on charging arcs do not exceed the limit of chargers available on a charging station. This constraint should ensure that, for example, Figure 2(c) should not occur if $u_2 = 1$ because technically both vehicle flows are concurrently using that charging station. To address this, the link flows $W_{igjh}$ are matched to enumerated path flows in constraint (11). There is one set of path flows for each charging station, so there are not many – for 4 charging levels there are 6 variables per charging station: e.g. {1-2, 1-3, 1-4, 2-3, 2-4, 3-4} where charge level 1 is lower than charge level 2. The path flows are used to ensure that path flow capacity is met in constraint (12). Lastly, $X_{igjh}$ and $Y_{jhm}$ are binary decision variables. Sayarshad and Chow (*7*) showed that $Y_{jhm}$ can generally be relaxed to a continuous variable between [0,1] since the piecewise linear constraint will generally be satisfied, which leads to a much more computationally efficient model. Arc flow $W_{igjh}$ is a non-negative integer decision variable of vehicle flow, and the path flows are continuous non-negative variables.

Note that in Eq. (5), $\rho_{\eta jm}$ is the coefficient of the utilization rate constraint, given a user-defined reliability rate $\eta$, m idle vehicles (servers) and b customers in a queue. The value of $\rho_{\eta jm}$ is obtained exogenously by solving the following Equation (17) (*7*), (*23-24*):

$$\sum_{k=0}^{m-1}((m-k)m!\,m^b/k!)\,(1/\rho^{m+b+1-k}) \geq 1/(1-\eta) \tag{17}$$

If the queueing constraint (Eq. (5)) is considered, the model represents a queueing constraint in a non-myopic case. Otherwise, the model is a myopic case without anticipating future queueing state in the system.

### 3.3 Solution method

The model simplifies to a p-median problem, which is known to be NP-complete (polynomial solvable if the "P" is exogenous, as noted in Garey and Johnson (*29*), which it is in this case). In this study we explore the differences in the new model formulation, so exact solution methods are needed. As such, we make use of commercial integer programming solvers to obtain a solution to the problem. For example, the intlinprog solver in MATLAB can obtain solutions using branch and bound or branch and cut algorithms.

## 4. NUMERICAL STUDY

Evaluation of the proposed model is conducted in two sets of replicable experiments. The first is with an illustrative example that serves to verify the model and demonstrate the capabilities to evaluate certain trade-offs. The second is with a set of 7 generated instances ranging in size from 10 to 1000 nodes, and 4 charging levels (up to 4,000 node-charges) to test the scalability of the model using commercial solvers.

### 4.1 Illustrative example

To illustrate the model, we consider a small network of 24 node-charges with 6 aligned nodes, extended up to 4 charging demand levels (20%, 40%, 60%, 80%) (see Figure 3). The goal is to test whether the proposed model can effectively rebalance the idle vehicles to meet all customer demand under available charging capacity constraints. Three idle vehicles with respective remaining charge levels are located at node-charge 3, 7, and 14. We consider 2 charging stations (node 2 and 6) with capacity $u$ per station. We test the model under three different capacities: $u = \{1,2,3\}$. Customer arrival rates are arbitrarily generated over all node-charges and fixed for the three scenarios. These rates are shown in the numbers over each node-charge (e.g. 3.8 customers/hr arrive at node-charge 15, which represents node 3 at charge level 60%+).

The model is implemented in MATLAB using a Dell Latitude E5470 laptop with win64 OS, Intel i5-6300U CPU, 2 Cores and 8GB memory. The Matlab mixed-integer linear programming solver (intlinprog) is used to solve the optimization problems for this study. The test instances are publicly available on the following data libraries: https://github.com/BUILTNYU.

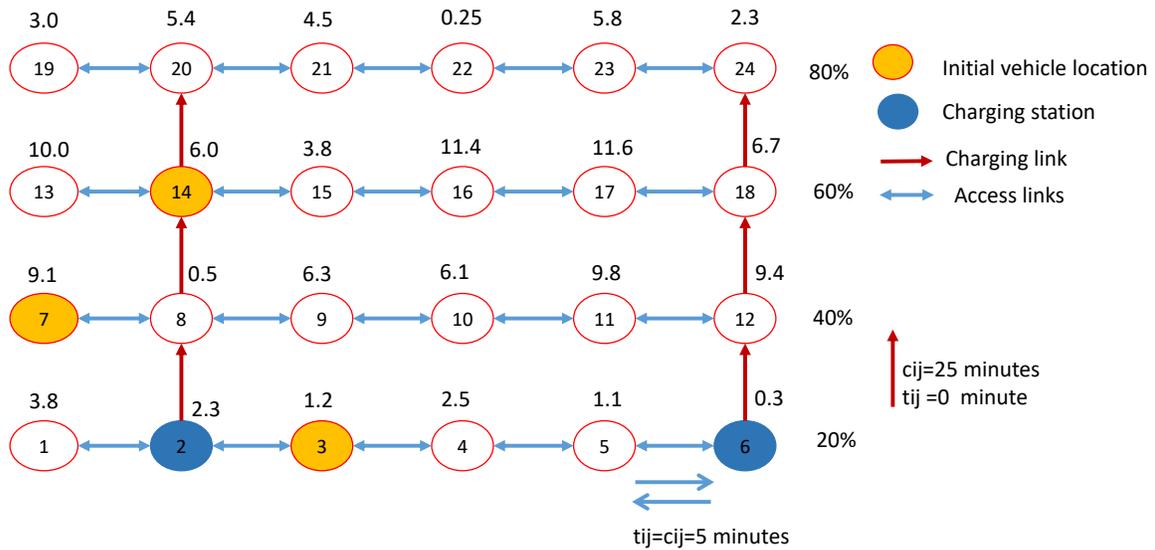

**FIGURE 3.** Illustration of a node-charge graph.

Figure 4 presents the computational result of rebalancing EV flows for the case of capacity $u = 3$. We see that all vehicles are rebalanced to the highest level 4 at different nodes which minimize total access cost of customers. The vehicles are assigned to use the nearest charging links to their destinations. The optimal objective value is $Z^* = 281.5$.

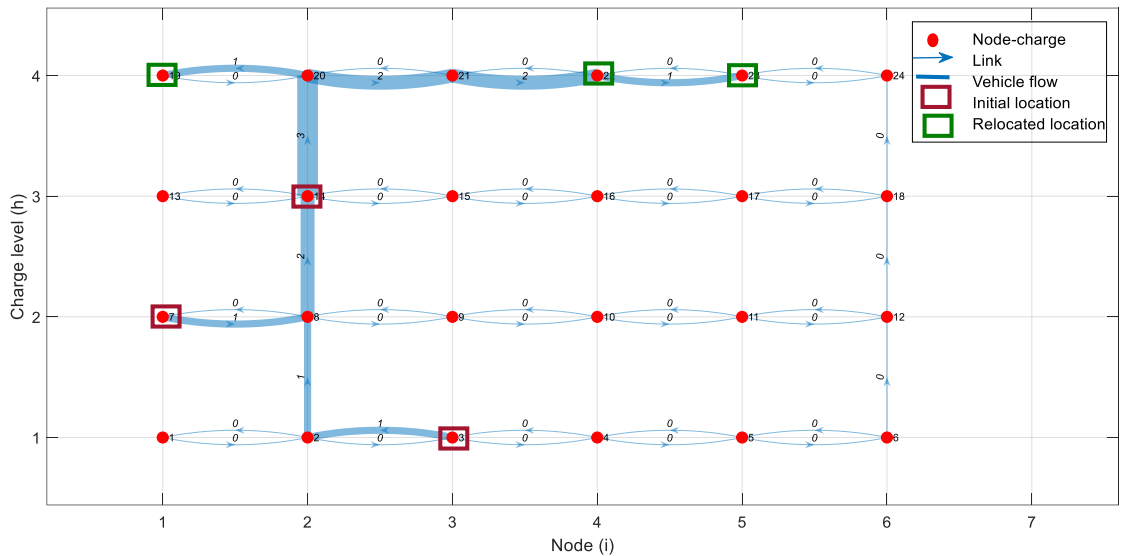

**FIGURE 4.** Rebalancing idle vehicle flows and its locations under charging station capacity = 3.

When reducing the capacity by 1 unit, we see the first charging station (i.e. node 2) is full capacitated by vehicles at node 7 and 14. The vehicle at node 3 moves farther to charge at node 6 and comes up to its destination at node 23 (see Figure 5(a)). The numbers over links represent the assigned flow on links. The obtained $Z^* = 281.5$ is the same as the preceding case.

When further reducing the capacity to $u = 1$ per station, there are only two chargers available in the system. The obtained solution shows the vehicles at node 7 and 3 coming up to the level 4 and the vehicle at node 14 is rebalanced on the same level to node 16. All charging capacity are used with least system cost (see Figure 5(b)). The obtained objective values is now $Z^* = 300.25$.

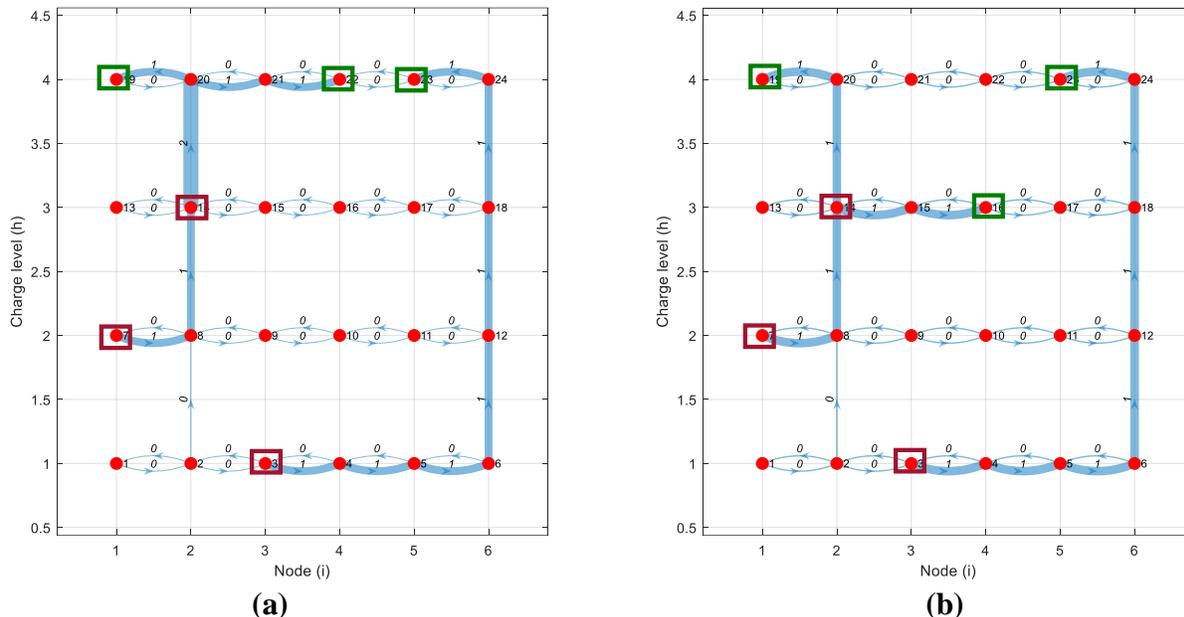

**FIGURE 5.** Rebalancing idle vehicle flows and its locations under charging station capacity constraints (Capacity =2 (a) and Capacity =1 (b)).

### 4.2 Large test cases with instances from 10 to 1000 nodes

To test the performance of the proposed model in large networks, we provide 7 test instances ranging from 40 node-charges to 4000 node-charges (4 charging demand levels: 20%, 40%, 60%, and 80%) (see Table 2). Customer demand at each node-charge is randomly generated from [0, 1]. The available idle vehicles are set as 10, and 4 charging stations with a capacity of 4 vehicles per station are considered. The parameter setting is shown in Table 1. The size of the test problems is shown in Table 2. We see the numbers of decision variables and of constraints (Eq. (2) - (11)) increase exponentially (162k decision variables for the 400 node-charge instance and 16 million constraints for the 4000 node-charge case).

First, we solve the myopic EV rebalancing problem by relaxing the queueing constraint (Eq. 5). The goal is to explore the computational performance of existing solver under a general personal computer. We see our model formulation allows finding optimal solution for 400 node-charges case in 10 seconds, around 40 minutes for 1600 node-charge case, and 2.14 hours for large network with 4000 node-charges. The result shows the performance of the model is quite fast and could be easily accelerated by using parallel computing on multithread computers for real-time operations.

Ma, Pantelidis, Chow                                                                                                               13Ma, Pantelidis, Chow                                                                                                               13

**TABLE 1 Reference parameter settings**

| Parameter | Value |
|---|---|
| $|N|$ | 10-1000 |
| $|H|$ | 4 |
| B and $|O|$ | 10 |
| $\lambda_{ig}$ | Random number drawn from [0,1] for each node-charge |
| $\theta$ | 0.2 |
| $C$ | 3 |
| $\rho_{jm}$ | (0.2236, 0.6416, 1.1576) |
| $|J|$ | 4 |
| $u_j$ | 4 |
| $M$ | 10000 |

**TABLE 2 Problem instance and computational times to solve myopic EV idle vehicle rebalancing problem**

| N | H | N × H | Num. of decision variables | Num. of equations | cputime (sec.) |
|---|---|---|---|---|---|
| 10 | 4 | 40 | 1 828 | 1 907 | <1 |
| 50 | 4 | 200 | 41 028 | 41 547 | 2 |
| 100 | 4 | 400 | 162 028 | 163 097 | 10 |
| 200 | 4 | 800 | 644 028 | 646 197 | 196 |
| 400 | 4 | 1600 | 2 568 028 | 2 572 397 | 2610 |
| 800 | 4 | 3200 | 10 256 028 | 10 264 797 | 7222 |
| 1000 | 4 | 4000 | 16 020 028 | 16 030 997 | 7724 |

For non-myopic case, the computational time depends on the complexity of Eq. (5), characterized by possible combinations of $Y_{jhm}$ and its relation between $\lambda$ and $\mu$. To illustrate this point, we explore the computational time on the case with 200 node-charges by varying the service rate $\mu_{jh}$ from $0.1N$ to $1.4N$ with an incremental $0.1N$ or $0.2N$ while keeping the same parameter setting. The result shows the computational time could be 2420 times higher when introducing the queueing constraint. Over the 12 tested $\mu_{jh}$, only three cases with $\mu_{jh} \geq N$ obtain optimal solution. No feasible solutions were found for the other cases. The result suggests further research is necessary in studying fast heuristics to find approximate solutions for large node-charge network in a non-myopic queueing constraint setting.

**TABLE 3 Influence of queueing constraints (eq. (5)) on computational time**

| Case | N × H | $\sum_{ig} \lambda_{ig}$ | $\mu_{jh}$ | With queueing constraint | |
|---|---|---|---|---|---|
| | | | | cputime (sec.) | Z* |
| 1 | 200 | 102.81 | 5 | 2 | NF |
| 2 | 200 | 102.81 | 10 | 2306 | NF |
| 3 | 200 | 102.81 | 15 | 2 | NF |
| 4 | 200 | 102.81 | 20 | 2 | NF |

| 5 | 200 | 102.81 | 25 | 6 | NF |
| --- | --- | --- | --- | --- | --- |
| 6 | 200 | 102.81 | 30 | 1816 | NF |
| 7 | 200 | 102.81 | 35 | 2600 | NF |
| 8 | 200 | 102.81 | 40 | 809 | NF |
| 9 | 200 | 102.81 | 45 | 4840 | NF |
| 10 | 200 | 102.81 | 50 | 22 | 831.206 |
| 11 | 200 | 102.81 | 60 | 9 | 829.305 |
| 12 | 200 | 102.81 | 70 | 3 | 827.305 |

Remark: NF means no feasible solutions. The reported cuptimes are the average of three runs for each case.

## 5. CONCLUSION

In this study, we propose a non-myopic idle vehicle rebalancing model in an electric carsharing system by considering customer stochastic charging demand and capacitated charging station constraints. To the best of our knowledge, this is the first facility relocation model formulation that considers queueing constraints applicable to EV charging. We formulate the problem as a p-median problem embedded with a capacitated minimum cost flow network problem on a node-charge graph to jointly determine the relocation and routing decisions of idle vehicles under available charging capacity.

The formulation on a two-dimensional node-charge graph allows us to explicitly consider a customer's charging demand profile and optimize rebalancing operations of idle vehicles accordingly. An illustrative example on a small network shows the assigned vehicle flow and rebalanced positions of idle vehicles minimize total designed objective while optimally utilizing available charging capacity in the system. We further test the performance of the model on a set of large test networks. The result shows our proposed mixed integer liner programming model formulation with commercial solvers allows us solving the large test instances efficiently with 10 seconds for 400 node-charges case (~162k decision variable and constraints each) and 2.14 hours for 4000 node-charges case (up to 16 million decision variables and constraints each) on a personal computer. This performance is promising given the size of the network it is run on (1000 nodes) and warrants the design of more efficient algorithms run on more efficient programming languages.

Several future extensions are currently ongoing. First, we are studying heuristics based on Lagrangian relaxation and greedy algorithms (e.g. (*30*)) approaches to find approximate solutions implying the non-myopic queueing constraint. Another research direction is studying the impact of the proposed model using realistic electric car sharing system data in New York City. The model will be evaluated in an online setting in comparison to benchmark methods.

## ACKNOWLEDGEMENTS

The work was supported by C2SMART University Transportation Center and the Luxembourg National Research Fund (INTER/MOBILITY/17/11588252).

**AUTHOR CONTRIBUTION STATEMENT**

The authors confirm contribution to the paper as follows: study conception and design, modeling, analysis, interpretation of results, draft manuscript preparation: T. Y. Ma, T. Pantelidis, J. Y. J. Chow. All authors reviewed the results and approved the final version of the manuscript.

**REFERENCES**

1. Martin, E., and S. Shaheen. Impacts of Car2go on Vehicle Ownership Modal Shift, Vehicle Miles Travelled, and Greenhouse Gas Emissions: An Analysis of Five North American Cities. Transportation Sustainability Research Center, UC Berkeley, California, 2016.
2. Bloomberg. Daimler, BMW reach a deal to merge car-sharing units. March 28, 2018.
3. Fan, W., R. Machemehl, and N. Lownes. Carsharing: Dynamic Decision-Making Problem for Vehicle Allocation. *Transportation Research Record: Journal of the Transportation Research Board*, 2008. 2063(1): 97–104.
4. Kek, A. G. H., Cheu, R. L., Meng, Q., and C. H. Fung. A Decision Support System for Vehicle Relocation Operations in Car Sharing Systems. *Transportation Research Part E*, 2009. 45(1): 149–158.
5. Nair, R., and E. Miller-Hooks. Fleet Management for Vehicle Sharing Operations. *Transportation Science*, 2011. 45(4): 451-566.
6. Nourinejad, M., S. Zhu, S. Bahrami, and M. J. Roorda. Vehicle Relocation and Staff Rebalancing in One-Way Carsharing Systems. *Transportation Research Part E*, 2015. 81: 98–113.
7. Sayarshad, H. R., and J. Y. J. Chow. Non-Myopic Relocation of Idle Mobility-on-Demand Vehicles as a Dynamic Location-Allocation-Queueing Problem. *Transportation Research Part E*, 2017. 106: 60–77.
8. Chow, J. Y. J., and H. R. Sayarshad. Symbiotic network design strategies in the presence of coexisting transportation networks. *Transportation Research Part B: Methodological*, 2014. 62: 13-34.
9. Jorge, D., G. Molnar, and G. H. de Almeida Correia. Trip pricing of one-way station-based carsharing networks with zone and time of day price variations. *Transportation Research Part B: Methodological*, 2015. 81: 461-482.
10. Waserhole, A., and V. Jost. Pricing in Vehicle Sharing Systems: Optimization in Queuing Networks with Product Forms. *EURO Journal on Transportation and Logistics*, 2016. 5(3): 293–320.
11. Clemente, M., M. P. Fanti, and W. Ukovich. Smart Management of Electric Vehicles Charging Operations: The Vehicle-to-Charging Station Assignment Problem. *IFAC Proceedings*, 2014. 19(3): 918–923.
12. Kaspi, M., T. Raviv, M. Tzur, and H. Galili. Regulating vehicle sharing systems through parking reservation policies: Analysis and performance bounds. *European Journal of Operational Research,* 2016. 251(3): 969-987.
13. Bruglieri, M., A. Colorni, and A. Luè. The relocation problem for the one-way electric vehicle sharing. *Networks*, 2014: 64(4), 292-305.
14. Hu, L., and Y. Liu. Joint design of parking capacities and fleet size for one-way station-based carsharing systems with road congestion constraints. *Transportation Research Part B: Methodological*, 2016. 93: 268-299.

15. Lufkin, B. 6 Electric Car-Sharing Programs Better Than a Billion Teslas on the Road. Gizmodo. https://gizmodo.com/6-electric-car-sharing-programs-better-than-a-billion-t-1766904615. March 31, 2016. Accessed Jul. 27, 2018.
16. Xu, M., Q. Meng, and Z. Liu. Electric vehicle fleet size and trip pricing for one-way carsharing services considering vehicle relocation and personnel assignment. *Transportation Research Part B: Methodological*, 2018. 111: 60-82.
17. Brandstätter, G., M. Kahr, and M. Leitner. Determining optimal locations for charging stations of electric car-sharing systems under stochastic demand. *Transportation Research Part B: Methodological*, 2017. 104: 17-35.
18. Boyacı, B., K. G. Zografos, and N. Geroliminis. An integrated optimization-simulation framework for vehicle and personnel relocations of electric carsharing systems with reservations. *Transportation Research Part B: Methodological*, 2017. 95: 214-237.
19. Li, X., J. Ma, J. Cui, A. Ghiasi, and F. Zhou. Design framework of large-scale one-way electric vehicle sharing systems: A continuum approximation model. *Transportation Research Part B: Methodological*, 2016. 88: 21-45.
20. Marianov, V., and D. Serra. Location–Allocation of Multiple-Server Service Centers with Constrained Queues or Waiting Times. *Annals of Operations Research*, 2002. 111: 35–50.
21. Barth, M., and M. Todd. Simulation model performance analysis of a multiple station shared vehicle system. *Transportation Research Part C: Emerging Technologies*, 1999. *7*(4): 237-259.
22. Febbraro, A., N. Sacco, and M. Saeednia. One-Way Carsharing: Solving the Relocation Problem. *Transportation Research Record: Journal of the Transportation Research Board,* 2012. 2319(1): 113–120.
23. Boyacı, B., K. G. Zografos, and N. Geroliminis. An Optimization Framework for the Development of Efficient One-Way Car Sharing Systems. *European Journal of Operational Research*, 2015. 240(3), 718–733.
24. Zhang, R., and M. Pavone (2016). Control of robotic mobility-on-demand systems: a queueing-theoretical perspective. *The International Journal of Robotics Research*, 35(1-3), 186-203.
25. Ma, T.Y., J.Y.J., Chow, and S. Rasulkhani, 2018. An integrated dynamic ridesharing dispatch and idle vehicle repositioning strategy on a bimodal transport network. *Proc. Transport Research Arena 2018*, Vienna, Austria.
26. Marianov, V., and C. ReVelleThe queueing maximal availability location problem: a model for the siting of emergency vehicles. *European Journal of Operational Research*, 1996. 93(1): 110-120.
27. Marianov, V., and D. Serra. Probabilistic, maximal covering location – allocation models for congested systems. *Journal of Regional Science,* 1998. 38: 401–424.
28. Jung, J., J. Y. J. Chow, R. Jayakrishnan, and J. Y. Park. Stochastic dynamic itinerary interception refueling location problem with queue delay for electric taxi charging stations. *Transportation Research Part C: Emerging Technologies*, 2014. 40: 123-142.
29. Garey, M. R., and D. S. Johnson. Computers and intractability: A guide to the theory of npcompleteness, ed. *Computers and Intractability*, 1979. 340.
30. Teitz, M. B., and P. Bart. Heuristic methods for estimating the generalized vertex median of a weighted graph. *Operations research*, 1968. 16(5): 955-961.